\documentclass[11pt]{article}
\usepackage{mathrsfs}
\usepackage[centertags]{amsmath}
\usepackage{amsfonts}
\usepackage{amssymb}
\usepackage{amsthm}
\usepackage{cases}
\usepackage{indentfirst}
\usepackage{epsfig}
\usepackage{graphicx}
\usepackage{color}
\usepackage{amsmath}

\date{}

\newtheorem{theorem}{Theorem}[section]

\newtheorem{lemma}{Lemma}[section]

\numberwithin{equation}{section}

\def\P{\operatorname*{\mathbf{P}}}

\def\R{\operatorname*{\mathbb{R}}}
\def\C{\operatorname*{\mathbb{C}}}

\def\SN{\operatorname*{\mathbb{SN}}}

\allowdisplaybreaks[4]

\begin{document}
\title{Uniform convergence rates of skew-normal extremes}
\author{\small{$^a$Qian Xiong\quad $^a$Zuoxiang Peng\quad $^b$Saralees Nadarajah}
\\
\\
\small{$^a$School of Mathematics and Statistics, Southwest University, 400715 Chongqing, China}
\\
\small{$^b$School of Mathematics, University of Manchester, Manchester, United Kingdom}}
\maketitle

\begin{quote}
{\bf Abstract.}~~Let $M_n=\max \left(X_1, X_2, \ldots, X_n \right)$
denote the partial maximum of an independent and identically distributed
skew-normal random sequence.
In this paper, the rate of uniform convergence of skew-normal extremes is derived.
It is shown that with optimal normalizing constants
the convergence rate of $\left(M_{n}-b_n\right)/a_n$
to its ultimate extreme value distribution is proportional to $1/\log n$.


\end{quote}

\section{Introduction}
\label{sec1}

Skew-normal distribution introduced by Azzalini (1985) is an
effective tool to model skewed data comparing to the normal distribution.
Applications of the skew-normal distribution
include areas such as climatology, biomedical sciences,
economics and finance.
Kim and Mallick (2004) presented a model based
on the skew-normal distribution for the prediction
of weekly rainfall in Korea.
Counsell et al. (2011) applied the skew-normal distribution
to deal with the data from a clinical psychiatry research environment.
Considering the effect of skewness and coskewness on asset valuation,
Carmichael and Co\"{e}n (2013) derived restrictions
imposed by the Euler equation of optimal portfolio diversification.
Zeller et al. (2016) established a mixture regression model by
assuming that the random errors follow a scale mixture of
skew-normal distributions.
With a skew-normal prior distribution for the spatial latent variables,
Hosseini et al. (2011) proposed approximate Bayesian methods
for inference and spatial prediction in a spatial generalized linear mixed model.
For more applications and case studies
involving skew-normal distributions, see Genton (2004).

Recently, probability properties  such as tail behavior
and asymptotics of skew-normal extremes have been studied.
Chang and Genton (2007) showed that  $F_{\lambda}$ belongs
to the domain of attraction of the Gumbel extreme
value distribution $\Lambda(x)=\exp (-\exp(-x) ), x\in\R$,
where $F_{\lambda}$ is the cumulative distribution
function (cdf) of the standard skew-normal random variable
with parameter $\lambda\in \R$ (written as $\SN(\lambda)$).
The  probability density function (pdf) of $\SN(\lambda)$ $f_{\lambda}$ is given by
\begin{eqnarray}
\label{eq1.1}
f_{\lambda}(x)= 2\phi(x)\Phi(\lambda x),
\quad
x\in\R,
\label{snpdf}
\end{eqnarray}
where $\phi (\cdot)$ and $\Phi(\cdot)$ denote, respectively, the pdf
and the cdf of a standard normal random variable.
For $\SN(0)$, the standard normal random variable, Mills' ratio and extreme
value distribution are known, cf. Leadbetter et al. (1983).
Higher-order expansions of the distribution and moments
of the extremes of $\SN (0)$ were studied by Nair (1981).
For $\lambda\neq0$, the Mills' ratios, the distributional tail
representation and the higher-order expansions of the extremes
of $\SN(\lambda)$ were studied by Liao et al. (2014b).
The higher-order expansions of moments of the extremes
of $\SN(\lambda)$, $\lambda \neq 0$ were studied by Liao et al. (2013a).
Liao et al. (2013b, 2014a) also considered the tail behaviors
and higher-order expansions of the distribution of extremes for
the log-skew-normal distribution.

The aim of this paper is to derive the uniform convergence
rates of skew-normal extremes for $\lambda\ne 0$.
For $\SN(0)$, the standard normal random variable,
Hall (1979) derived the optimal uniform convergence
rate of $\Phi^n \left(\widetilde{a}_nx+\widetilde{b}_n\right)$ to $\Lambda(x)$, i.e.,
\begin{eqnarray}
\label{eq1.2}
\frac {\C_1}{\log n}<\sup_{x\in\R} \left|
\Phi^n \left(\widetilde{a}_nx+\widetilde{b}_n\right) - \Lambda(x) \right |<\frac {\C_2}{\log n}
\end{eqnarray}
for some positive constants $\C_1$ and $\C_2$,
where the normalizing constant $\widetilde{b}_n$ is the solution of  equation
\begin{eqnarray*}
2\pi \widetilde{b}^2_n \exp\left(\widetilde{b}_n^2\right)=n^2
\end{eqnarray*}
and $\widetilde{a}_n=\widetilde{b}_n^{-1}$.
For $\SN(0)$, Gasull et al. (2015a) gave more effective
normalizing constants $a_{n}$ and $b_{n}$ through the Lambert W function.
Gasull et al. (2015b) illustrates another
application of the Lambert W function to decide
on normalizing constants for gamma and other Weibull-like distributions.
For other work related to the convergence rates of
distributions of normalized order statistics,
see Liao and Peng (2012) for the log-normal distribution,
Peng et al. (2010) and  Vasudeva et al. (2014) for the general error distribution.

In order to derive the uniform convergence rates of skew-normal extremes,
we chose the optimal normalizing constants according
to the sign of $\lambda$ with $\lambda\ne 0$.
For $\lambda>0$, let $b_n$ be the solution of the equation
\begin{eqnarray}
\label{eq1.3}
\sqrt{\frac {\pi}{2}}b_n \exp\left({\frac {1}{2}b_n^2}\right)=n
\end{eqnarray}
and
\begin{eqnarray}
\label{eq1.4}
a_n=b_n^{-1}.
\end{eqnarray}
For$\lambda<0$, let $b_n>0$ be the solution of the equation
\begin{eqnarray}
\label{eq1.5}
\pi|\lambda|\left(1+\lambda^2\right)
b_n^2 \exp\left({\frac {\left(1+\lambda^2\right)b_n^2}{2}}\right)=n
\end{eqnarray}
and
\begin{eqnarray}
\label{eq1.6}
a_n=\left(1+\lambda^2\right)^{-1}b_n^{-1}.
\end{eqnarray}

The rest of this paper is organized as follows.
Section \ref{sec2} gives the main results.
Some auxiliary lemmas and all proofs are presented in Section \ref{sec3}.

\section{Main results}
\label{sec2}

In this section, we provide the main results.
Theorem \ref{th2.1} shows that the limit distribution of
normalized maxima for the skew-normal distribution is
the Gumbel extreme value distribution $\Lambda(x)$.
Theorem \ref{th2.2} and Theorem \ref{th2.3} determine the
rates of uniform convergence of skew-normal extremes.
Note that the choice of normalizing constants are determined
according to the sign of $\lambda$.

\begin{theorem}
\label{th2.1}
Let $M_n$ denote the partial maximum of
independent and identical $\SN(\lambda)$ random variables
with the pdf $f_{\lambda}$ given by \eqref{eq1.1}.
Then,
\begin{eqnarray}
\label{eq2.1}
\P\left(M_n\leq a_n x+b_n\right)\rightarrow\Lambda(x),
\quad
x\in \R
\end{eqnarray}
as $n\rightarrow\infty$, where the normalizing
constants $b_n$ and $a_n$ are given by \eqref{eq1.3} and \eqref{eq1.4}
for $\lambda>0$, and by \eqref{eq1.5} and \eqref{eq1.6} for $\lambda<0$.
\end{theorem}

\begin{theorem}
\label{th2.2}
For $\lambda>0$, there exist positive constants $\C$ and $\C_\lambda$,
independent of $n$, such that for all $n\geq9$,
\begin{eqnarray}\label{eq2.2}
\frac {\C}{\log n}<\sup_{x\in\R} \left | F_\lambda^n \left(a_nx+b_n\right)-\Lambda(x) \right |<\frac {\C_\lambda}{\log n},
\end{eqnarray}
where $b_n$ and $a_n$ are given by \eqref{eq1.3} and \eqref{eq1.4}, respectively.
\end{theorem}

\begin{theorem}
\label{th2.3}
For $\lambda<0$, there exist positive constants $\C_{\lambda}^{'}$ and $\C_{\lambda}^{''}$, independent of $n$, such that for all $n\geq n_0(\lambda)$,
\begin{eqnarray}
\label{eq2.3}
\frac {\C_{\lambda}^{'}}{\log n}<\sup_{x\in\R}
\left |F_\lambda^n \left(a_nx+b_n\right)-\Lambda(x) \right |<
\frac {\C_{\lambda}^{''}}{\log n},
\end{eqnarray}
where $b_n$ and $a_n$ are given by \eqref{eq1.5} and \eqref{eq1.6}, respectively.
\end{theorem}

\section{Proofs}
\label{sec3}

In order to prove the main results, we first give some auxiliary lemmas.
The first one is about the distributional tail representation
of the skew-normal distribution  due to Lemma 3.1 in Xiong and Peng (2018).
The remaining lemmas provide  inequalities on distributional tails
of the normal and skew-normal distributions.

\begin{lemma}
\label{lem3.1}
Let $F_{\lambda}$ and $f_{\lambda}$ denote
the cdf and pdf of the $\SN(\lambda)$ distribution.
For large $x$, we have
\begin{itemize}

\item[(i)]
for $\lambda>0$,
\begin{eqnarray}
\label{eq3.1}
1-F_{\lambda}(x)=\frac {2\phi(x)}{x}\left[1-x^{-2}+O\left(x^{-4}\right) \right].
\end{eqnarray}

\item[(ii)]
for $\lambda <0$,
\begin{eqnarray}
\label{eq3.2}
1-F_{\lambda}(x)=\frac {e^{-\frac {\left(1+\lambda^2\right)x^2}{2}}}{-\pi\lambda \left(1+\lambda^2\right)x^2}
\left[1-\frac {1+3\lambda^2}{\lambda^2 \left(1+\lambda^2\right)}x^{-2}+
O\left(x^{-4}\right) \right].
\end{eqnarray}
\end{itemize}
\end{lemma}

\begin{lemma}
\label{lem3.2}
Let $\phi (x)$ and $\Phi(x)$ denote respectively the pdf and the cdf of a standard normal random variable.
For all $x>0$, we have
\begin{eqnarray}
\label{eq3.3}
\frac {\phi(x)}{x} \left(1-x^{-2}\right) <1-\Phi(x)<\frac {\phi(x)}{x}.
\end{eqnarray}
\end{lemma}

\noindent
{\bf Proof. }
The proof is straightforward by integration by parts, see equations (6)-(9) in Hall (1979).

\begin{lemma}
\label{lem3.3}
Let $F_\lambda$ denote the cdf of $\SN(\lambda)$ and
let $\phi (x)$ denote the pdf of a standard normal random variable.
For all $x>0$, we have
\begin{itemize}

\item[(i)]
for $\lambda>0$,
\begin{eqnarray}
\label{eq3.4}
\frac {2\phi(x)}{x}\left[ 1-\left(1+\frac {1}{\lambda^2\sqrt{2\pi e}}\right)x^{-2}\right]<1-F_\lambda(x)<\frac {2\phi(x)}{x};
\end{eqnarray}

\item[(ii)]
for $\lambda<0$,
\begin{eqnarray}
\label{eq3.5}
\frac {2\phi(x)\phi(\lambda x)}{|\lambda| \left(1+\lambda^2\right)x^2}
\left[ 1-\frac {\left(1+\lambda^2\right)^2}{\lambda^2}x^{-2}\right] <
1-F_\lambda(x)
<\frac {2\phi(x)\phi(\lambda x)}{|\lambda|\left(1+\lambda^2\right)x^2}.
\end{eqnarray}
\end{itemize}
\end{lemma}

\noindent
{\bf Proof. }
In the case of $\lambda>0$, for any $x>0$ we have
\begin{eqnarray}
\label{eq3.6}
1-F_\lambda(x)=\int_x^\infty 2\phi(x)\Phi(\lambda t)dt
<\int_x^\infty 2\phi(x)dt=2\left[ 1-\Phi(x)\right] <\frac {2\phi(x)}{x}.
\end{eqnarray}
By integration by parts and Lemma \ref{lem3.2}, we have
\begin{eqnarray}
\label{eq3.7}
1-F_\lambda(x)
&=&
\frac {f_\lambda(x)}{x}-\int_x^\infty 2\phi(t)\Phi(\lambda t)t^{-2}dt+
\frac {\lambda}{\pi}\int_x^\infty t^{-1}e^{-\frac {\left(1+\lambda^2\right)t^2}{2}}dt
\nonumber
\\
&>&
\frac {f_\lambda(x)}{x}-\int_x^\infty 2\phi(t)\Phi(\lambda t)t^{-2}dt
\nonumber
\\
&>&
\frac {f_\lambda(x)}{x}-\int_x^\infty 2\phi(t)t^{-2}dt
\nonumber
\\
&>&
\frac {f_\lambda(x)}{x}-x^{-2}\int_x^\infty 2\phi(t)dt
\nonumber
\\
&=&
\frac {2\phi(x)\Phi(\lambda x)}{x}-2x^{-2}\left[1-\Phi(x)\right]
\nonumber
\\
&>&
\frac {2\phi(x)}{x}\left[1-\frac {\phi(\lambda x)}{\lambda x}\right] -
2x^{-2}\cdot\frac {\phi(x)}{x}\nonumber
\\
&=&
\frac {2\phi(x)}{x}\left\{ 1-\left[ 1+\frac {x\phi(\lambda x)}{\lambda}\right] x^{-2}\right\}
\nonumber
\\
&>&
\frac {2\phi(x)}{x}\left[1-\left(1+\frac {1}{\lambda^2\sqrt{2\pi e}}\right)x^{-2}\right].
\end{eqnarray}
The last inequality was obtained by bounding the
function $x\exp{\left(-\frac {\lambda^2x^2}{2}\right)}$.
Combining \eqref{eq3.6} with \eqref{eq3.7}, we can derive \eqref{eq3.4}.

In the case of $\lambda<0$, for any $x>0$,
by Lemma \ref{lem3.2} we have
\begin{eqnarray}
\label{eq3.8}
1-F_\lambda(x)
&=&
\int_x^\infty 2\phi(t)\Phi(\lambda t)dt
\nonumber
\\
&=&
\int_x^\infty 2\phi(t)\left[ 1-\Phi(|\lambda|t)\right] dt
\nonumber
\\
&<&
\int_x^\infty 2\phi(x)\cdot\frac {\phi(|\lambda|t)}{|\lambda|t}dt
\nonumber
\\
&=&
\frac {1}{\pi|\lambda|}\int_x^\infty t^{-1}
e^{-\frac {\left(1+\lambda^2\right)t^2}{2}}dt
\nonumber
\\
&=&
\frac {2\phi(x)\phi(\lambda x)}{|\lambda|\left(1+\lambda^2\right)x^2}-
\frac {2}{\pi|\lambda|\left(1+\lambda^2\right)}
\int_x^\infty t^{-3}e^{-\frac {\left(1+\lambda^2\right)t^2}{2}}dt
\nonumber
\\
&<&
\frac {2\phi(x)\phi(\lambda x)}{|\lambda| \left(1+\lambda^2\right)x^2}.
\end{eqnarray}
By integration by parts and Lemma \ref{lem3.2}, we have
\begin{eqnarray}
\label{eq3.9}
1-F_\lambda(x)
&=&
\frac {f_\lambda(x)}{x}-\frac {f_\lambda(x)}{x^3}
+6\int_x^\infty \phi(t)\Phi(\lambda t)t^{-4}dt+
\frac {|\lambda|}{\pi}\int_x^\infty t^{-3}e^{-\frac {\left(1+\lambda^2\right)t^2}{2}}dt
\nonumber
\\
&&
-\frac {f_\lambda(x)}{x}\cdot\frac {|\lambda|}{1+\lambda^2}\cdot
\frac {\phi(\lambda x)}{\Phi(\lambda x)}x^{-1}
+\frac {2|\lambda|}{\left(1+\lambda^2\right)\pi}\int_x^\infty t^{-3}
e^{-\frac {\left(1+\lambda^2\right)t^2}{2}}dt
\nonumber
\\
&>&
\frac {f_\lambda(x)}{x}-\frac {f_\lambda(x)}{x^3}-
\frac {f_\lambda(x)}{x}\cdot\frac {|\lambda|}{1+\lambda^2}
\cdot\frac {\phi(\lambda x)}{\Phi(\lambda x)}x^{-1}
\nonumber
\\
&=&
\frac {2\phi(x)}{x}\left[\Phi(\lambda x) \left(1-x^{-2}\right) -
\frac {|\lambda|}{1+\lambda^2}\cdot\frac {\phi(\lambda x)}{x}\right]
\nonumber
\\
&>&
\frac {2\phi(x)}{x}\left[ \frac {\phi(\lambda x)}{|\lambda|x}\left(1-\lambda^{-2}x^{-2}\right)
\left(1-x^{-2}\right) - \frac {|\lambda|}{1+\lambda^2}\cdot
\frac {\phi(\lambda x)}{x}\right]
\nonumber
\\
&>&
\frac {2\phi(x)\phi(\lambda x)}{|\lambda|x^2}
\left[ \left(1-\lambda^{-2}x^{-2}\right) \left(1-x^{-2}\right) -
\frac {\lambda^2}{1+\lambda^2}\right]
\nonumber
\\
&>&
\frac {2\phi(x)\phi(\lambda x)}{|\lambda| \left(1+\lambda^2\right)x^2}
\left[ 1-\frac {\left(1+\lambda^2\right)^2}{\lambda^2}x^{-2}\right].
\end{eqnarray}
Combining \eqref{eq3.8} with \eqref{eq3.9}, we can derive \eqref{eq3.5}.
\qed

\noindent
{\bf Proof of Theorem \ref{th2.1}.}
We first consider the case of $\lambda>0$.
If $n$ is sufficiently large then $a_nx+b_n>0$,
with $b_n$ and $a_n$ satisfying \eqref{eq1.3} and \eqref{eq1.4},
and so by Lemma \ref{lem3.1} we have
\begin{eqnarray*}
n\left[1-F_\lambda \left(a_nx+b_n\right) \right]
&\sim&
n\cdot\frac {2\phi \left(a_nx+b_n\right)}{a_nx+b_n}
\\
&=&
\left(1+a_n^2x\right)^{-1}e^{-\frac {1}{2}a_n^2x^2}\cdot e^{-x}
\\
&\to&
e^{-x}
\end{eqnarray*}
as $n\to\infty$.
Then the result follows from Theorem 1.5.1 of Leadbetter et al. (1983).

Similarly, for $\lambda<0$,
if $n$ is sufficiently large then $a_nx+b_n>0$
with $b_n$ and $a_n$ satisfying \eqref{eq1.5} and \eqref{eq1.6},
and so by Lemma \ref{lem3.1} we have
\begin{eqnarray*}
n\left[1-F_\lambda \left(a_nx+b_n\right) \right]
&\sim&
n\cdot\frac {e^{-\frac {\left(1+\lambda^2\right)
\left(a_nx+b_n\right)^2}{2}}}{\pi|\lambda|\left(1+\lambda^2\right)\left(a_nx+b_n\right)^2}
\\
&=&
\left[ 1+\left(1+\lambda^2\right)a_n^2x\right]^{-2}
e^{-\frac {\left(1+\lambda^2\right)a_n^2x^2}{2}}\cdot e^{-x}
\\
&\to&
e^{-x}
\end{eqnarray*}
as $n\to\infty$.
Then the result follows from Theorem 1.5.1 of Leadbetter et al. (1983).
\qed

\noindent
{\bf Proof of Theorem \ref{th2.2}.}
First note that for $\lambda>0$, sufficiently large $n$ implies
$a_nx+b_n>0$ with $b_n$ and $a_n$ satisfying \eqref{eq1.3} and \eqref{eq1.4}.
Then for large $n$ and $k\in\R$, we have
\begin{eqnarray}
\label{eq3.10}
\left(a_nx+b_n\right)^k =
b_n^k \left(1+a_n^2x\right)^k=b_n^k\left[ 1+ka_n^2x+O\left(a_n^4\right)\right].
\end{eqnarray}
Applying \eqref{eq3.10} to $k=-1$, $k=2$ and $k=-2$, we obtain
\begin{eqnarray}
\label{eq3.11}
\frac {2\phi \left(a_nx+b_n\right)}{a_nx+b_n} =
n^{-1}e^{-x}\left[ 1-a_n^2\left(x+\frac {1}{2}x^2\right)+O\left(a_n^4\right)\right]
\end{eqnarray}
and
\begin{eqnarray}
\label{eq3.12}
1 - \left(a_nx+b_n\right)^{-2} = 1-b_n^{-2}
\left[1-2a_n^2x+O\left(a_n^4\right)\right] = 1-a_n^2+O\left(a_n^4\right).
\end{eqnarray}
Combining \eqref{eq3.11} with \eqref{eq3.12}, we have
\begin{eqnarray}
\label{eq3.13}
\frac {2\phi \left(a_nx+b_n\right)}{a_nx+b_n}
\left[ 1 - \left(a_nx+b_n\right)^{-2}\right] =
n^{-1}e^{-x}\left[1-a_n^2\left(1+x+\frac {1}{2}x^2\right)+O\left(a_n^4\right)\right].
\end{eqnarray}
Therefore by Lemma \ref{lem3.1} and the
fact that $\log(1-x)=-x \left[ 1+O(x) \right]$ as $x\to0$, we have
\begin{eqnarray*}
F^n_\lambda \left(a_nx+b_n\right)-\Lambda(x)
&=&
\left\{1-\frac {2\phi\left(a_nx+b_n\right)}{a_nx+b_n}
\left[1-\left(a_nx+b_n\right)^{-2}+
O\left(\left(a_nx+b_n\right)^{-4}\right)\right]\right\}^n-\Lambda(x)
\\
&=&
\left\{1-n^{-1}e^{-x}\left[1-a_n^2\left(1+x+\frac {1}{2}x^2\right)
+O\left(a_n^4\right)\right]\right\}^n-\Lambda(x)
\\
&=&
\exp\left\{n\log\left\{1-n^{-1}e^{-x}\left[1-a_n^2\left(1+x+\frac {1}{2}x^2\right)+
O\left(a_n^4\right)\right]\right\}\right\}-\Lambda(x)
\\
&=&
\Lambda(x)\exp\left\{e^{-x}\left[a_n^2\left(1+x+\frac {1}{2}x^2\right)+
O\left(a_n^4\right)\right]\right\}-\Lambda(x)
\\
&=&
\Lambda(x)\left\{1+e^{-x}\left[a_n^2\left(1+x+\frac {1}{2}x^2\right)+
O\left(a_n^4\right)\right] \right\} -\Lambda(x)
\\
&=&
\Lambda(x)e^{-x}\left[a_n^2\left(1+x+\frac {1}{2}x^2\right)+O\left(a_n^4\right)\right].
\end{eqnarray*}
Further, by \eqref{eq1.3}, we have
\begin{eqnarray}
\label{eq3.14}
\log\frac {\pi}{2}+2\log b_n+b_n^2=2\log n.
\end{eqnarray}
It follows at once that $b_n^2\sim2\log n$.
Noting that $a_n=b_n^{-1}$, then we can
obtain the left-hand-side inequality in \eqref{eq2.2}.
So, the rest is to show that for all $n\geq9$,
there exists a positive constant $\C_\lambda$ such that
\begin{eqnarray*}
\sup_{x\in\R} \left |F_\lambda^n \left( a_nx+b_n \right)-\Lambda(x) \right |<\frac {\C_\lambda}{\log n}.
\end{eqnarray*}

For $n\geq2$, \eqref{eq3.14} implies that
\begin{eqnarray}
\label{eq3.15}
b_n^2<2\log n,
\end{eqnarray}
so that
\begin{eqnarray}
\label{eq3.16}
2\log b_n<\log2 +\log\log n.
\end{eqnarray}
Combining \eqref{eq3.16} with \eqref{eq3.14}, we obtain
\begin{eqnarray*}
b_n^2>2\log n-\log\pi-\log\log n,
\end{eqnarray*}
and hence for $n\geq9$, we obtain
\begin{eqnarray}
\label{eq3.17}
\frac {b_n^2}{\log n}
&>&
2-\frac {\log\pi}{\log n}-\frac {\log\log n}{\log n}
\nonumber
\\
&>&
2-\frac {\log\pi}{\log 9}-\frac {1}{e}
\nonumber
\\
&>&
1.1.
\end{eqnarray}
The second inequality is obtained by bounding
the function $(\log x)^{-1}(\log\log x)$.
Since $a_n=b_n^{-1}$, the inequality \eqref{eq3.17}
implies that for $n\geq9$, $a_n^2<\frac {1}{1.1\log n}$,
and so it suffices to prove that for $n\geq9$,
\begin{eqnarray}
\label{eq3.18}
\sup_{x\in\R} \left |F_\lambda^n \left(a_nx+b_n\right)-\Lambda(x) \right |<{\mathbb{C}_\lambda}a_n^2.
\end{eqnarray}
We will prove the result by showing
\begin{eqnarray}
\label{eq3.19}
\sup_{0\leq x<\infty}
\left | F_\lambda^n \left(a_nx+b_n\right)-\Lambda(x) \right |<{\mathbb{C}_{1,\lambda}}a_n^2,
\end{eqnarray}
\begin{eqnarray}
\label{eq3.20}
\sup_{-c_n<x<0}
\left | F_\lambda^n \left(a_nx+b_n\right) - \Lambda(x) \right |<
{\mathbb{C}_{2,\lambda}}a_n^2
\end{eqnarray}
and
\begin{eqnarray}
\label{eq3.21}
\sup_{-\infty<x\leq -c_n}
\left| F_\lambda^n \left(a_nx+b_n\right)-\Lambda(x) \right |<
{\mathbb{C}_{3,\lambda}}a_n^2,
\end{eqnarray}
where $c_n=\log\log b_n^2>0$ for $n\geq9$.

The following bounds are needed for the rest of the proof:
\begin{eqnarray}
\label{eq3.22}
1.69<b_9<1.70,
\end{eqnarray}
\begin{eqnarray}
\label{eq3.23}
\sup_{n\geq9} \left(1-a_n^2c_n\right)^{-1}<1.11,
\end{eqnarray}
\begin{eqnarray}
\label{eq3.24}
\sup_{n\geq9} a_n^2\log b_n^2<0.37,
\end{eqnarray}
\begin{eqnarray}
\label{eq3.25}
\sup_{n\geq9} a_n^2 \left(\log b_n^2\right)^2<0.55,
\end{eqnarray}
\begin{eqnarray}
\label{eq3.26}
\sup_{n\geq9} n^{-1}\log b_n^2<0.17
\end{eqnarray}
and
\begin{eqnarray}
\label{eq3.27}
\sup_{n\geq9} b_n^3e^{-\frac {1}{2}b_n^2}<1.16.
\end{eqnarray}
The inequality \eqref{eq3.22} follows from \eqref{eq1.3}
and \eqref{eq3.26} follows from \eqref{eq3.15}.
Inequalities \eqref{eq3.23}, \eqref{eq3.24}, \eqref{eq3.25}
and \eqref{eq3.27} are obtained by bounding the functions
$x^{-1}\log\log x$, $x^{-1}\log x$, $x^{-1}(\log x)^2$ and
$x^3e^{-\frac {1}{2}x^2}$, respectively.

Let $\Psi_{n,\lambda}(x)=1-F_\lambda \left(a_nx+b_n\right)$, then
\begin{eqnarray}\label{eq3.28}
n\log F_\lambda \left(a_nx+b_n\right)=
n\log\left[{1-\Psi_{n,\lambda}(x)}\right]=-n\Psi_{n,\lambda}(x)-R_{n,\lambda}(x),
\end{eqnarray}
where
\begin{eqnarray}\label{eq3.29}
0<R_{n,\lambda}(x)<
\frac {n\Psi^2_{n,\lambda}(x)}{2\left[1-\Psi_{n,\lambda}(x)\right]}.
\end{eqnarray}
If $x>-c_n$, by \eqref{eq1.3}, \eqref{eq3.23}, \eqref{eq3.26} and Lemma \ref{lem3.3}, we have
\begin{eqnarray}
\label{eq3.30}
\Psi_{n,\lambda}(x)
&<&
\Psi_{n,\lambda}\left(-c_n\right)
\nonumber
\\
&=&
1-F_\lambda\left(b_n-a_nc_n\right)
\nonumber
\\
&<&
\frac {2\phi \left(b_n-a_nc_n\right)}{b_n-a_nc_n}
\nonumber
\\
&=&
\sqrt{\frac {2}{\pi}}b_n^{-1}
\left(1-a_n^2c_n\right)^{-1}e^{-\frac {1}{2}b_n^2\left(1-a_n^2c_n\right)^2}
\nonumber
\\
&=&
n^{-1}\left(1-a_n^2c_n\right)^{-1}
\left(\log b_n^2\right) e^{-\frac {1}{2}a_n^2c_n^2}
\nonumber
\\
&<&
\left(1-a_n^2c_n\right)^{-1}
\left( n^{-1}\log b_n^2\right)
\nonumber
\\
&<&
0.1887.
\end{eqnarray}
From \eqref{eq3.29}, \eqref{eq3.30} \eqref{eq1.3} and \eqref{eq1.4}, we can see that
\begin{eqnarray}
\label{eq3.31}
R_{n,\lambda}(x)
&<&
\frac {n\left[\left(1-a_n^2c_n\right)^{-1}
\left(n^{-1}\log b_n^2\right)\right]^2}{2(1-0.1887)}
\nonumber
\\
&=&
\frac {n^{-1}\left(1-a_n^2c_n\right)^{-2}\left(\log b_n^2\right)^2}{1.6226}
\nonumber
\\
&=&
\frac {\sqrt{\frac {2}{\pi}}b_n^{-1}
e^{-\frac {1}{2}b_n^2}
\left(1-a_n^2c_n\right)^{-2}\left(\log b_n^2\right)^2}{1.6226}
\nonumber
\\
&=&
\frac {\left[\sqrt{\frac {2}{\pi}}\left(1-a_n^2c_n\right)^{-2}\right]
\left(b_n^3e^{-\frac {1}{2}b_n^2}\right)
\left[a_n^2\left(\log b_n^2\right)^2\right] a_n^2}{1.6226}
\nonumber
\\
&<&
0.39a_n^2.
\end{eqnarray}
The last inequality follows from \eqref{eq3.23}, \eqref{eq3.25},
\eqref{eq3.27}, and hence for $n\geq9$, we have
\begin{eqnarray}
\label{eq3.32}
\left |e^{-R_{n,\lambda}(x)}-1\right |=1-e^{-R_{n,\lambda}(x)}<R_{n,\lambda}(x)<0.39a_n^2.
\end{eqnarray}

Let $A_{n,\lambda}(x)=e^{-n\Psi_{n,\lambda}(x)+e^{-x}}$
and $B_{n,\lambda}(x)=e^{-R_{n,\lambda}(x)}$.
The inequality \eqref{eq3.32} implies that
\begin{eqnarray}
\label{eq3.33}
\left | F^n_\lambda \left(a_nx+b_n\right)-\Lambda(x) \right |
&<&
\Lambda(x) \left |A_{n,\lambda}(x)-1 \right | + \left | B_{n,\lambda}(x)-1 \right |
\nonumber
\\
&<&
\Lambda(x) \left |A_{n,\lambda}(x)-1 \right |+0.39a_n^2
\end{eqnarray}
if $x>-c_n$.

We first show that \eqref{eq3.19} holds.
Note that $0<A_{n,\lambda}(x)\to 1$ as $x\to \infty$ and
\begin{eqnarray*}
A^{'}_{n,\lambda}(x) =
A_{n,\lambda}(x)e^{-x}\left[e^{-\frac {1}{2}a_n^2x^2}
\Phi\left(\lambda \left( a_nx+b_n \right) \right)-1\right]\le 0
\end{eqnarray*}
for $x\geq0$.
Hence, it follows from \eqref{eq3.22}, \eqref{eq1.4} and Lemma \ref{lem3.3} that
\begin{eqnarray}
\label{eq3.34}
\sup_{x\geq0} \left |A_{n,\lambda}(x)-1 \right |
&=&
A_{n,\lambda}(0)-1
\nonumber
\\
&<&
\exp\left\{-n\frac {2\phi \left(b_n\right)}{b_n}
\left[1-\left(1+\frac {1}{\lambda^2\sqrt{2\pi e}}\right)
b_n^{-2}\right]+1\right\}-1
\nonumber
\\
&<&
\left(1+\frac {1}{\lambda^2\sqrt{2\pi e}}\right)b_n^{-2}\exp\left\{\left(1+\frac {1}{\lambda^2\sqrt{2\pi e}}\right)b_n^{-2}\right\}
\nonumber
\\
&<&
\left(1+\frac {1}{\lambda^2\sqrt{2\pi e}}\right)
\exp\left\{1.69^{-2}\left(1+\frac {1}{\lambda^2\sqrt{2\pi e}}\right)\right\} a_n^2.
\end{eqnarray}
The second inequality is obtained by observing that $e^x<1+xe^x$ for $x>0$,
and the last inequality is obtained by the monotonicity of $b_n$ (i.e. $b_n\geq b_9$ for $n\geq9$).
Combining \eqref{eq3.33} with \eqref{eq3.34}, we finish the proof of \eqref{eq3.19}.

Next we prove that \eqref{eq3.20} holds as $-c_n<x<0$.
By \eqref{eq1.4}, \eqref{eq3.23} and the fact that $e^x>1+x$ for $x \in \R$, we have
\begin{eqnarray}\label{eq3.35}
&&
-e^{-\frac {1}{2}a_n^2x^2}\left[1-\left(1+\frac {1}{\lambda^2\sqrt{2\pi e}}\right)
\left(a_nx+b_n\right)^{-2}\right]+1+a_n^2x
\nonumber
\\
&<&
-\left(1-\frac {1}{2}a_n^2x^2\right)
\left[1-\left(1+\frac {1}{\lambda^2\sqrt{2\pi e}}\right)
\left(a_nx+b_n\right)^{-2}\right]+1+a_n^2x
\nonumber
\\
&=&
a_n^2\left(1+a_n^2x\right)^{-2}
\left(1+\frac {1}{\lambda^2\sqrt{2\pi e}}\right)+
\frac {1}{2}a_n^2x^2-\frac {1}{2}a_n^4x^2
\left(1+a_n^2x\right)^{-2}
\left(1+\frac {1}{\lambda^2\sqrt{2\pi e}}\right)+a_n^2x
\nonumber
\\
&<&
a_n^2\left(1+a_n^2x\right)^{-2}
\left(1+\frac {1}{\lambda^2\sqrt{2\pi e}}\right)+
\frac {1}{2}a_n^2x^2
\nonumber
\\
&<&
a_n^2\left[ \left(1-a_n^2c_n\right)^{-2}
\left(1+\frac {1}{\lambda^2\sqrt{2\pi e}}\right)+\frac {1}{2}x^2\right]
\nonumber
\\
&<&
a_n^2\left[1.11^2\left(1+\frac {1}{\lambda^2\sqrt{2\pi e}}\right)+
\frac {1}{2}x^2\right].
\end{eqnarray}
Let $h_{n,\lambda}(x)=-n\Psi_{n,\lambda}(x)+e^{-x}$.
From Lemma \ref{lem3.3}, \eqref{eq1.3}, \eqref{eq1.4} and \eqref{eq3.35}, we have
\begin{eqnarray}
\label{eq3.36}
h_{n,\lambda}(x)
&<&
-n\frac {2\phi \left(a_nx+b_n\right)}{a_nx+b_n}
\left[1-\left(1+\frac {1}{\lambda^2\sqrt{2\pi e}}\right)
\left(a_nx+b_n\right)^{-2}\right]+e^{-x}
\nonumber
\\
&=&
-n\sqrt{\frac {2}{\pi}}b_n^{-1}
\left(1+a_n^2x\right)^{-1} e^{-\frac {1}{2}\left(a_nx+b_n\right)^2}
\left[1-\left(1+\frac {1}{\lambda^2\sqrt{2\pi e}}\right)
\left(a_nx+b_n\right)^{-2}\right]+e^{-x}
\nonumber
\\
&=&
\left(1+a_n^2x\right)^{-1}e^{-x}\left\{-e^{-\frac {1}{2}a_n^2x^2}
\left[1-\left(1+\frac {1}{\lambda^2\sqrt{2\pi e}}\right)
\left(a_nx+b_n\right)^{-2}\right]+1+a_n^2x\right\}
\nonumber
\\
&<&
\left(1+a_n^2x\right)^{-1}e^{-x}a_n^2
\left[1.11^2\left(1+\frac {1}{\lambda^2\sqrt{2\pi e}}\right)+\frac {1}{2}x^2\right].
\end{eqnarray}
Further by Lemma \ref{lem3.3}, \eqref{eq1.3} and \eqref{eq1.4}, we have
\begin{eqnarray}
\label{eq3.37}
h_{n,\lambda}(x)
&>&
-n\frac {2\phi\left(a_nx+b_n\right)}{a_nx+b_n}+e^{-x}
\nonumber
\\
&=&
-n\sqrt{\frac {2}{\pi}}b_n^{-1}
\left(1+a_n^2x\right)^{-1} e^{-\frac {1}{2}\left(a_nx+b_n\right)^2}+e^{-x}
\nonumber
\\
&=&
\left(1+a_n^2x\right)^{-1}e^{-x}\left(-e^{-\frac {1}{2}a_n^2x^2}+1+a_n^2x\right)
\nonumber
\\
&>&
\left(1+a_n^2x\right)^{-1}e^{-x}a_n^2|x|.
\end{eqnarray}
Hence, for $-c_n<x<0$, it follows from \eqref{eq3.36} and \eqref{eq3.37} that
\begin{eqnarray}
\label{eq3.38}
\left |h_{n,\lambda}(x) \right |
&<&
\left(1+a_n^2x\right)^{-1}e^{-x}a_n^2
\left[1.11^2\left(1+\frac {1}{\lambda^2\sqrt{2\pi e}}\right)+
\frac {1}{2}x^2+|x|\right]
\nonumber
\\
&<&
\left(1-a_n^2c_n\right)^{-1}e^{c_n}
a_n^2\left[1.11^2\left(1+\frac {1}{\lambda^2\sqrt{2\pi e}}\right)+
\frac {1}{2}c_n^2+c_n\right]
\nonumber
\\
&<&
\left(1-a_n^2c_n\right)^{-1}
\left[1.11^2\left(1+\frac {1}{\lambda^2\sqrt{2\pi e}}\right)
\left(a_n^2 \log b_n^2 \right)+\frac {3}{2}a_n^2 \left(\log b_n^2 \right)^2\right]
\nonumber
\\
&<&
\mathbb{C}_{2,\lambda}^{'}.
\end{eqnarray}
The third inequality holds by the facts $\log x>(\log\log x)^2$
for $x>e$ and $\log x>\log\log x$ for $x>1$,
and the last inequality holds by \eqref{eq3.23}-\eqref{eq3.25}.
Noting that $\left |e^x-1 \right |<|x|e^{|x|}$ for $x\in\R$
and $e^x>1+x+\frac {1}{2}x^2$ for $x>0$, we have
\begin{eqnarray}
\label{eq3.39}
\Lambda(x) \left |A_{n,\lambda(x)}-1 \right |
&=&
\Lambda(x) \left |e^{h_{n,\lambda}(x)}-1 \right |
\nonumber
\\
&<&
\Lambda(x) \left | h_{n,\lambda}(x) \right |
e^{\left |h_{n,\lambda}(x) \right |}
\nonumber
\\
&<&
\Lambda(x) \left(1+a_n^2x\right)^{-1}
e^{-x}a_n^2\left[1.11^2\left(1+\frac {1}{\lambda^2\sqrt{2\pi e}}\right)+
\frac {1}{2}x^2+|x|\right]e^{\mathbb{C}_{2,\lambda}^{'}}
\nonumber
\\
&<&
\Lambda(x)\left(1-a_n^2c_n\right)^{-1}
e^{-x}a_n^2\left[1.11^2\left(1+\frac {1}{\lambda^2\sqrt{2\pi e}}\right)+
\frac {1}{2}x^2+|x|\right]e^{\mathbb{C}_{2,\lambda}^{'}}
\nonumber
\\
&=&
a_n^2\left(1-a_n^2c_n\right)^{-1}
\left[1.11^2\left(1+\frac {1}{\lambda^2\sqrt{2\pi e}}\right)+
\frac {1}{2}x^2+|x|\right]e^{{-e^{-x}-x+\mathbb{C}_{2,\lambda}^{'}}}
\nonumber
\\
&<&
a_n^2\left(1-a_n^2c_n\right)^{-1}\left[1.11^2
\left(1+\frac {1}{\lambda^2\sqrt{2\pi e}}\right)+
\frac {1}{2}x^2+|x|\right]e^{{-\frac {1}{2}x^2+\mathbb{C}_{2,\lambda}^{'}}-1}
\nonumber
\\
&=&
a_n^2\left(1-a_n^2c_n\right)^{-1}\left[1.11^2
\left(1+\frac {1}{\lambda^2\sqrt{2\pi e}}\right)
e^{-\frac {1}{2}x^2}+\frac {1}{2}x^2
e^{-\frac {1}{2}x^2}+|x|e^{-\frac {1}{2}x^2}\right]
e^{{\mathbb{C}_{2,\lambda}^{'}}-1}
\nonumber
\\
&<&
1.11 a_n^2 \left[1.11^2\left(1+\frac {1}{\lambda^2\sqrt{2\pi e}}\right)+
1+1\right]e^{{\mathbb{C}_{2,\lambda}^{'}}-1}.
\end{eqnarray}
Putting \eqref{eq3.39} into \eqref{eq3.33}, we can establish \eqref{eq3.20}.

The last step is to show that \eqref{eq3.21} holds.
For $-\infty<x\leq-c_n$,  \eqref{eq1.4} implies
\begin{eqnarray}
\label{eq3.40}
0\leq\Lambda(x)\leq\Lambda \left(-c_n\right)=a_n^2.
\end{eqnarray}
Since $e^x>1+x$ for $x \in \R$, we have
\begin{eqnarray}\label{eq3.41}
&&
\left(1-a_n^2c_n\right)^{-1} e^{-\frac {1}{2}a_n^2c_n^2}
\left[1-a_n^2\left(1-a_n^2c_n\right)^{-2}
\left(1+\frac {1}{\lambda^2\sqrt{2\pi e}}\right)\right]
\nonumber
\\
&>&
e^{-\frac {1}{2}a_n^2c_n^2}\left[1-a_n^2\left(1-a_n^2c_n\right)^{-2}
\left(1+\frac {1}{\lambda^2\sqrt{2\pi e}}\right)\right]
\nonumber
\\
&>&
\left(1-\frac {1}{2}a_n^2c_n^2\right)
\left[1-a_n^2\left(1-a_n^2c_n\right)^{-2}
\left(1+\frac {1}{\lambda^2\sqrt{2\pi e}}\right)\right]
\nonumber
\\
&>&
1-a_n^2\left[\frac {1}{2}c_n^2+\left(1-a_n^2c_n\right)^{-2}
\left(1+\frac {1}{\lambda^2\sqrt{2\pi e}}\right)\right].
\end{eqnarray}
Thus, by \eqref{eq1.3}, \eqref{eq1.4}, \eqref{eq3.23},
\eqref{eq3.24}, \eqref{eq3.25}, \eqref{eq3.41} and Lemma \ref{lem3.3}, we have
\begin{eqnarray}
\label{eq3.42}
F_\lambda^n\left(a_nx+b_n\right)
&\leq&
F_\lambda^n \left(b_n-a_nc_n\right)
\nonumber
\\
&=&
\left\{ 1-\left[ 1 - F_\lambda \left(b_n-a_nc_n\right) \right] \right\}^n
\nonumber
\\
&<&
\left\{1-\frac {2\phi\left(b_n-a_nc_n\right)}{b_n-a_nc_n}
\left[1-\left(1+\frac {1}{\lambda^2\sqrt{2\pi e}}\right)
\left(b_n-a_nc_n\right)^{-2}\right] \right\}^n
\nonumber
\\
&=&
\left\{1-\sqrt{\frac {2}{\pi}}b_n^{-1}
\left(1-a_n^2c_n\right)^{-1}
e^{-\frac {1}{2}\left(b_n-a_nc_n\right)^2}
\left[1-a_n^2\left(1-a_n^2c_n\right)^{-2}
\left(1+\frac {1}{\lambda^2\sqrt{2\pi e}}\right)\right]\right\}^n
\nonumber
\\
&=&
\left\{1-n^{-1}e^{c_n}
\left(1-a_n^2c_n\right)^{-1}
e^{-\frac {1}{2}a_n^2c_n^2}\left[ 1-a_n^2\left(1-a_n^2c_n\right)^{-2}
\left(1+\frac {1}{\lambda^2\sqrt{2\pi e}}\right)\right] \right\}^n
\nonumber
\\
&<&
\exp\left\{-e^{c_n}\left(1-a_n^2c_n\right)^{-1}
e^{-\frac {1}{2}a_n^2c_n^2}
\left[ 1-a_n^2\left(1-a_n^2c_n\right)^{-2}
\left(1+\frac {1}{\lambda^2\sqrt{2\pi e}}\right)\right]\right\}
\nonumber
\\
&<&
\exp\left\{-e^{c_n}\left[1-a_n^2\left(\frac {1}{2}c_n^2+
\left(1-a_n^2c_n\right)^{-2}\left(1+\frac {1}{\lambda^2\sqrt{2\pi e}}\right)
\right)\right]\right\}
\nonumber
\\
&<&
a_n^2\exp\left\{\frac {1}{2}a_n^2\left(\log b_n^2\right)^2+
\left[a_n^2 \left(\log b_n^2\right)\right]
\left(1-a_n^2c_n\right)^{-2}
\left(1+\frac {1}{\lambda^2\sqrt{2\pi e}}\right)
\right\}
\nonumber
\\
&<&
\mathbb{C}^{'}_{3,\lambda}a_n^2
\end{eqnarray}
due to $\left(1-\frac {x}{n}\right)^n<e^{-x}$ for $x\leq n$ and
$\log x>(\log\log x)^2$ for $x>e$.
Combining \eqref{eq3.40} with \eqref{eq3.42}, we see that \eqref{eq3.21} holds.
The proof is complete.
\qed

\noindent
{\bf Proof of Theorem \ref{th2.3}.}
For $\lambda<0$, if $n$ is sufficiently large then $a_nx+b_n>0$
with $b_n$ and $a_n$ satisfying \eqref{eq1.5} and \eqref{eq1.6}.
Writing $z_n=a_nx+b_n$, for large $n$ and $k\in\R$, we have
\begin{eqnarray}\label{eq3.43}
z_n^k=\left(a_nx+b_n\right)^k=b_n^k
\left[1+\left(1+\lambda^2\right)a_n^2x\right]^k=
b_n^k\left[1+k\left(1+\lambda^2\right)a_n^2x+O\left(a_n^4\right)\right].
\end{eqnarray}
Applying \eqref{eq3.43} with $k=2$ and $k=-2$, we have
\begin{eqnarray}
\label{eq3.44}
\frac {e^{-\frac {\left(1+\lambda^2\right)z_n^2}{2}}}
{-\pi\lambda \left(1+\lambda^2\right)z_n^2} =
n^{-1}e^{-x}\left[1 - \left(1+\lambda^2\right)
a_n^2\left(2x+\frac {1}{2}x^2\right)+O \left(a_n^4\right)\right]
\end{eqnarray}
and
\begin{eqnarray}
\label{eq3.45}
1-\frac {1+3\lambda^2}{\lambda^2 \left(1+\lambda^2\right)}z_n^{-2} =
1-\frac {\left(1+3\lambda^2\right) \left(1+\lambda^2\right)}{\lambda^2}
a_n^2+O \left(a_n^4\right).
\end{eqnarray}
Combining \eqref{eq3.44} with \eqref{eq3.45}, we have
\begin{eqnarray}
\label{eq3.46}
\frac {e^{-\frac {\left(1+\lambda^2\right)z_n^2}{2}}}
{-\pi\lambda \left(1+\lambda^2\right)z_n^2}
\left[1-\frac {1+3\lambda^2}{\lambda^2\left(1+\lambda^2\right)}z_n^{-2}\right] =
n^{-1}e^{-x}\left[1 - \left(1+\lambda^2\right)
a_n^2\left(\frac {1+3\lambda^2}{\lambda^2}+
2x+\frac {1}{2}x^2\right)+O\left(a_n^4\right)\right].
\end{eqnarray}
Therefore, by Lemma \ref{lem3.1}, we have
\begin{eqnarray*}
F_\lambda^n \left(a_nx+b_n\right) -\Lambda(x)
&=&
\left\{1-\frac {e^{-\frac {\left(1+\lambda^2\right)z_n^2}{2}}}{-\pi\lambda \left(1+\lambda^2\right)z_n^2}
\left[1-\frac {1+3\lambda^2}{\lambda^2 \left(1+\lambda^2\right)}z_n^{-2}+
O \left(z_n^{-4}\right)\right]\right\}^n-\Lambda(x)
\nonumber
\\
&=&
\left\{1-n^{-1}e^{-x}\left[1-\left(1+\lambda^2\right)
a_n^2\left(\frac {1+3\lambda^2}{\lambda^2}+2x+
\frac {1}{2}x^2\right)+O \left(a_n^4\right)\right]\right\}^n-\Lambda(x)
\nonumber
\\
&=&
\Lambda(x)e^{-x}\left[\left(1+\lambda^2\right)
a_n^2\left(\frac {1+3\lambda^2}{\lambda^2}+
2x+\frac {1}{2}x^2\right)+O\left(a_n^4\right)\right].
\end{eqnarray*}
Further, by \eqref{eq1.5} we have
\begin{eqnarray}
\label{eq3.47}
\log\left[\pi|\lambda|\left(1+\lambda^2\right)\right]+
2\log b_n+\frac {\left(1+\lambda^2\right)b_n^2}{2}=\log n,
\end{eqnarray}
implying that $\left(1+\lambda^2\right)b_n^2\sim2\log n$.
Noting that $a_n=\left(1+\lambda^2\right)^{-1}b_n^{-1}$,
we can obtain the left-hand-side inequality in \eqref{eq2.3}.
The rest is to show that for all $n>n_0(\lambda)$,
there exists a positive constant $\C_\lambda^{''}$ such that
\begin{eqnarray*}
\sup_{x\in\R}\left | F_\lambda \left(a_nx+b_n\right)^n-\Lambda(x) \right | <
\frac {\mathbb{C}_\lambda^{''}}{\log n}.
\end{eqnarray*}

For $n>n_0(\lambda)$,  \eqref{eq3.47} implies
\begin{eqnarray}
\label{eq3.48}
b_n^2<\frac {2}{\left(1+\lambda^2\right)}\log n,
\end{eqnarray}
so that
\begin{eqnarray}
\label{eq3.49}
2\log b_n<\log\frac {2}{\left(1+\lambda^2\right)}+\log\log n.
\end{eqnarray}
Combining \eqref{eq3.49} with \eqref{eq3.47}, we have
\begin{eqnarray*}
\left(1+\lambda^2\right)b_n^2>2\log n-2\log \left(2\pi|\lambda|\right) - 2\log\log n
\end{eqnarray*}
and
\begin{eqnarray}
\label{eq3.50}
\frac {\left(1+\lambda^2\right)b_n^2}{\log n}
&>&
2-\frac {2\log \left(2\pi|\lambda|\right)}{\log n}-\frac {2\log\log n}{\log n}
\nonumber
\\
&>&
2-\frac {2\log \left(2\pi|\lambda|\right)}{\log n_0(\lambda)}-\frac {2}{e}
\nonumber
\\
&=&
c_0,
\end{eqnarray}
where $c_0$ is a positive constant and the last inequality is obtained
by bounding the function $(\log x)^{-1}(\log\log x)$.
Since $a_n = \left(1+\lambda^2\right)^{-1}b_n^{-1}$,
\eqref{eq3.50} implies that for $n\geq n_0(\lambda)$,
$a_n^2<\frac {1}{c_0 \left(1+\lambda^2\right)\log n}$,
and so it suffices to prove that for $n\geq n_0(\lambda)$
\begin{eqnarray}
\label{eq3.51}
\sup_{x\in\R}
\left | F_\lambda^n \left(a_nx+b_n\right) - \Lambda(x) \right |<{\mathbb{C}_\lambda^{''}}a_n^2.
\end{eqnarray}
We will now turn to prove the following inequalities:
\begin{eqnarray}
\label{eq3.52}
\sup_{0\leq x<\infty}
\left |F_\lambda^n \left(a_nx+b_n\right) - \Lambda(x) \right | < {\mathbb{C}_{1,\lambda}^{''}}a_n^2,
\end{eqnarray}
\begin{eqnarray}
\label{eq3.53}
\sup_{-d_n<x<0}
\left |F_\lambda^n \left(a_nx+b_n\right) - \Lambda(x) \right | <{\mathbb{C}_{2,\lambda}^{''}}a_n^2
\end{eqnarray}
and
\begin{eqnarray}
\label{eq3.54}
\sup_{-\infty<x\leq -d_n}
\left | F_\lambda^n \left(a_nx+b_n\right) - \Lambda(x) \right| <{\mathbb{C}_{3,\lambda}^{''}}a_n^2,
\end{eqnarray}
where $d_n=\log\log \left[ \left(1+\lambda^2\right)b_n^2\right]$ and $d_n>0$ for $n\geq n_0(\lambda)$.

The following bounds are needed in our proof:
\begin{eqnarray}
\label{eq3.55}
c_1<b_{n_0(\lambda)},
\end{eqnarray}
\begin{eqnarray}
\label{eq3.56}
\sup_{n\geq2}
\left[ 1 - \left(1+\lambda^2\right)a_n^2d_n\right]^{-1}<1.11,
\end{eqnarray}
\begin{eqnarray}
\label{eq3.57}
\sup_{n\geq2}
\left(1+\lambda^2\right)a_n^2\log \left[ \left(1+\lambda^2\right)b_n^2\right] <0.37,
\end{eqnarray}
\begin{eqnarray}
\label{eq3.58}
\sup_{n\geq n_0(\lambda)}
\left(1+\lambda^2\right)
a_n^2\left\{ \log \left[ \left(1+\lambda^2\right) b_n^2\right] \right\}^2<0.55,
\end{eqnarray}
\begin{eqnarray}
\label{eq3.59}
\sup_{n\geq2}
n^{-1}\log \left[ \left(1+\lambda^2\right)b_n^2\right] <0.27
\end{eqnarray}
and
\begin{eqnarray}
\label{eq3.60}
\sup_{n\geq2} \left(1+\lambda^2\right)
b_n^2e^{-\frac {\left(1+\lambda^2\right)b_n^2}{2}}<0.74,
\end{eqnarray}
where $c_1$ is a positive constant and
the inequality \eqref{eq3.59} follows from \eqref{eq3.48},
and \eqref{eq3.56}, \eqref{eq3.57}, \eqref{eq3.58}
and \eqref{eq3.60} are obtained by bounding the
functions $x^{-1}\log\log x$, $x^{-1}\log x$, $x^{-1}(\log x)^2$
and $xe^{-\frac {1}{2}x}$, respectively.

If $x>-d_n$, by \eqref{eq1.5}, \eqref{eq3.29},
\eqref{eq3.56}, \eqref{eq3.59} and Lemma \ref{lem3.3}, we have
\begin{eqnarray}
\label{eq3.63}
\Psi_{n,\lambda}(x)
&<&
\Psi_{n,\lambda}\left(-d_n\right)
\nonumber
\\
&<&
\frac {2\phi\left(b_n-a_nd_n\right)\phi\left(\lambda \left(b_n-a_nd_n\right)\right)}
{|\lambda| \left(1+\lambda^2\right) \left(b_n-a_nd_n\right)^2}
\nonumber
\\
&=&
n^{-1}\left[ 1-\left(1+\lambda^2\right)a_n^2d_n\right]^{-2}
\left\{ \log\left[ \left(1+\lambda^2\right)b_n^2\right] \right\}
e^{-\frac {\left(1+\lambda^2\right)a_n^2d_n^2}{2}}
\nonumber
\\
&<&
\left[ 1 - \left(1+\lambda^2\right)a_n^2d_n\right]^{-2}
\left\{ n^{-1}\log\left[ \left(1+\lambda^2\right)b_n^2\right] \right\}
\nonumber
\\
&<&
0.332667
\end{eqnarray}
and
\begin{eqnarray}
\label{eq3.64}
R_{n,\lambda}(x)
&<&
\frac {n\left\{ \left[1-\left(1+\lambda^2\right)a_n^2d_n\right]^{-2}
\left\{n^{-1}\log\left[\left(1+\lambda^2\right)b_n^2\right]\right\} \right\}^2}
{2(1-0.332667)}
\nonumber
\\
&<&
\frac {1}{\pi|\lambda|}
\left\{\left[1 - \left(1+\lambda^2\right)a_n^2d_n\right]^{-4}\right\}
\left[\left(1+\lambda^2\right)a_n^2
\left\{ \log\left[ \left(1+\lambda^2\right)b_n^2\right]\right\}^2\right]
\nonumber
\\
&&
\cdot
\left[\left(1+\lambda^2\right)b_n^2
e^{-\frac {\left(1+\lambda^2\right)b_n^2}{2}}\right]
\left(1+\lambda^2\right)a_n^2
\nonumber
\\
&<&
\mathbb{C}_4^{''}a_n^2,
\end{eqnarray}
where $\Psi_{n,\lambda}(x)$ and $R_{n,\lambda}(x)$
are given by \eqref{eq3.28},
i.e., $\Psi_{n,\lambda}(x) =1-F_\lambda \left(a_nx+b_n\right)$ and
$R_{n,\lambda}(x)=-n\log F_\lambda\left(a_nx+b_n\right)-n\Psi_{n,\lambda}(x)$
with $b_{n}$ and $a_{n}$ given by \eqref{eq1.5} and \eqref{eq1.6}, respectively.
Hence, for $n\geq n_0(\lambda)$, we obtain
\begin{eqnarray}\label{eq3.65}
\left |e^{-R_{n,\lambda}(x)}-1 \right | =
1-e^{-R_{n,\lambda}(x)}<R_{n,\lambda}(x)<\mathbb{C}_4^{''}a_n^2
\end{eqnarray}
by using the inequality $e^x>1+x$ for $x \in \R$.

Let $A_{n,\lambda}(x)=e^{-n\Psi_{n,\lambda}(x)+e^{-x}}$
and $B_{n,\lambda}(x)=e^{-R_{n,\lambda}(x)}$.
\eqref{eq3.65} implies that
\begin{eqnarray}
\label{eq3.66}
\left | F^n_\lambda \left(a_nx+b_n\right) - \Lambda(x) \right |
&<&
\Lambda(x) \left | A_{n,\lambda}(x)-1 \right | + \left | B_{n,\lambda}(x)-1 \right |
\nonumber
\\
&<&
\Lambda(x) \left | A_{n,\lambda}(x)-1 \right |+\mathbb{C}_4^{''}a_n^2
\end{eqnarray}
as $x>-d_n$.

Now we prove \eqref{eq3.52}-\eqref{eq3.54} in turn.
To prove \eqref{eq3.52}, noting that $0<A_{n,\lambda}(x)\to 1$
as $x\to \infty$ and by Lemma \ref{lem3.2}, we have
\begin{eqnarray*}
A^{'}_{n,\lambda}(x)
&=&
A_{n,\lambda}(x)\left[-e^{-x}+na_n f_\lambda \left(a_nx+b_n\right)\right]
\\
&<&
A_{n,\lambda}(x)\left[-e^{-x}+na_n 2\phi\left(a_nx+b_n\right)
\frac {\phi\left(|\lambda|\left(a_nx+b_n\right)\right)}{|\lambda|\left(a_nx+b_n\right)}\right]
\\
&=&
A_{n,\lambda}(x)e^{-x}
\left\{e^{-\frac {1}{2}\left(1+\lambda^2\right)a_n^2x^2}
\left[1+\left(1+\lambda^2\right)a_n^2x\right]^{-1}-1\right\}
\\
&\leq&
0
\end{eqnarray*}
for $x\geq0$.
Hence, by \eqref{eq3.55}, \eqref{eq1.4} and Lemma \ref{lem3.3}, we have
\begin{eqnarray}\label{eq3.67}
\sup_{x\geq0}
\left |A_{n,\lambda}(x)-1\right |
&=&
A_{n,\lambda}(0)-1
\nonumber
\\
&<&
\exp\left\{-n\frac {2\phi \left(b_n\right)
\phi \left(\lambda b_n\right)}
{|\lambda| \left(1+\lambda^2\right) b_n}
\left[1-\frac {\left(1+\lambda^2\right)^2}{\lambda^2}b_n^{-2}\right]+1\right\}-1
\nonumber
\\
&=&
\exp\left\{\frac {\left(1+\lambda^2\right)^2}{\lambda^2}b_n^{-2}\right\}-1
\nonumber
\\
&<&
\frac {\left(1+\lambda^2\right)^2}{\lambda^2}b_n^{-2}
\exp\left\{\frac {\left(1+\lambda^2\right)^2}{\lambda^2}b_n^{-2}\right\}
\nonumber
\\
&<&
\frac {\left(1+\lambda^2\right)^4}{\lambda^2}
\exp\left\{\frac {\left(1+\lambda^2\right)^2}{\lambda^2} c_1^{-2}\right\}a_n^2
\end{eqnarray}
due to $e^x<1+xe^x$ for $x>0$ and the monotonicity of $b_n$
(i.e. $b_n\geq b_{n_0(\lambda)}$ for $n\geq n_0(\lambda)$).
Combining \eqref{eq3.66} with \eqref{eq3.67}, we prove \eqref{eq3.52}.

Before proving \eqref{eq3.53}, we need the following inequalities.
By \eqref{eq1.5} and the fact that $e^x>1+x$ for $x \in \R$, we have
\begin{eqnarray}
\label{eq3.68}
&&
-e^{-\frac {1}{2}\left(1+\lambda^2\right)
a_n^2x^2}\left[1-\frac {\left(1+\lambda^2\right)^2}{\lambda^2}
\left(a_nx+b_n\right)^{-2}\right]+\left[ 1+\left(1+\lambda^2\right)a_n^2x\right]^2
\nonumber
\\
&<&
-\left[ 1-\frac {1}{2}\left(1+\lambda^2\right)a_n^2x^2\right]
\left[1-\frac {\left(1+\lambda^2\right)^2}{\lambda^2}
\left(a_nx+b_n\right)^{-2}\right]+
\left[1+\left(1+\lambda^2\right)a_n^2x\right]^2
\nonumber
\\
&<&
\left(1+\lambda^2\right)a_n^2\left\{\frac {\left(1+\lambda^2\right)^3}{\lambda^2}
\left[1+\left(1+\lambda^2\right)a_n^2x\right]^{-2}+
\frac {1}{2}x^2+\left(1+\lambda^2\right)a_n^2x^2\right\}.
\end{eqnarray}
Let $h_{n,\lambda}(x)=-n\Psi_{n,\lambda}(x)+e^{-x}$.
By  \eqref{eq1.5}, \eqref{eq1.6}, \eqref{eq3.68} and Lemma \ref{lem3.3}, we have
\begin{eqnarray}
\label{eq3.69}
h_{n,\lambda}(x)
&<&
-n\frac {2\phi \left(a_nx+b_n\right)\phi\left(\lambda \left(a_nx+b_n\right)\right)}{|\lambda|
\left(1+\lambda^2\right)\left(a_nx+b_n\right)^2}
\left[1-\frac {\left(1+\lambda^2\right)^2}{\lambda^2}\left(a_nx+b_n\right)^{-2}\right]+e^{-x}
\nonumber
\\
&<&
\left[1+\left(1+\lambda^2\right)a_n^2x\right]^{-2}e^{-x}
\left(1+\lambda^2\right)a_n^2
\Bigg\{\frac {\left(1+\lambda^2\right)^3}{\lambda^2}
\left[ 1+\left(1+\lambda^2\right)a_n^2x\right]^{-2}
\nonumber
\\
&&
+\frac {1}{2}x^2+\left(1+\lambda^2\right)a_n^2x^2\Bigg\}
\end{eqnarray}
and
\begin{eqnarray}
\label{eq3.70}
h_{n,\lambda}(x)
&>&
-n\frac {2\phi\left(a_nx+b_n\right)
\phi\left(\lambda \left(a_nx+b_n\right)\right)}
{|\lambda| \left(1+\lambda^2\right) \left(a_nx+b_n\right)^2}+e^{-x}
\nonumber
\\
&>&
\left[ 1+\left(1+\lambda^2\right)a_n^2x\right]^{-2}
e^{-x}\left( 1+\lambda^2 \right)
a_n^2\left[2x+\left(1+\lambda^2\right)a_n^2x^2\right].
\end{eqnarray}
Hence, for $-d_n<x<0$, we note from \eqref{eq3.69} and \eqref{eq3.70}  that
{\small\begin{eqnarray}
\label{eq3.71}
&&
\left |h_{n,\lambda}(x) \right |
\nonumber
\\
&<&
\left[1+ \left(1+\lambda^2\right)a_n^2x\right]^{-2}
e^{-x}\left(1+\lambda^2\right)
a_n^2\left\{\frac {\left(1+\lambda^2\right)^3}{\lambda^2}
\left[ 1+\left(1+\lambda^2\right)a_n^2x\right]^{-2}+
\frac {1}{2}x^2+2\left(1+\lambda^2\right)a_n^2x^2+2|x|\right\}
\nonumber
\\
&<&
\left[1-\left(1+\lambda^2\right)a_n^2d_n\right]^{-2}
e^{d_n} \left(1+\lambda^2\right)
a_n^2\left\{ \frac {\left(1+\lambda^2\right)^3}{\lambda^2}
\left[1 - \left(1+\lambda^2\right)a_n^2d_n\right]^{-2}
+\frac {1}{2}d_n^2+2\left(1+\lambda^2\right)a_n^2d_n^2+2d_n\right\}
\nonumber
\\
&=&
\left[1 - \left(1+\lambda^2\right)a_n^2d_n\right]^{-2}
\Bigg[ \frac {\left(1+\lambda^2\right)^3}{\lambda^2}
\left[1 - \left(1+\lambda^2\right)a_n^2d_n\right]^{-2}
\left\{ \left(1+\lambda^2\right)a_n^2
\log\left[ \left(1+\lambda^2\right) b_n^2\right] \right\}
\nonumber
\\
&&
+\frac {1}{2}\left\{ \left(1+\lambda^2\right)
a_n^2\log\left[ \left(1+\lambda^2\right)b_n^2\right] \right\}
\left\{ \log\log\left[ \left(1+\lambda^2\right)b_n^2\right] \right\}^2
\nonumber
\\
&&
+2\left\{ \left(1+\lambda^2\right)^2a_n^4
\log\left[ \left(1+\lambda^2\right) b_n^2\right] \right\}
\left\{ \log\log\left[ \left(1+\lambda^2\right) b_n^2\right] \right\}^2
\nonumber
\\
&&
+2\left\{ \left(1+\lambda^2\right)a_n^2\log\left[ \left(1+\lambda^2\right)b_n^2\right]
\right\} \left\{ \log\log\left[ \left(1+\lambda^2\right)b_n^2\right] \right\} \Bigg]
\nonumber
\\
&<&
\mathbb{C}_{5,\lambda}^{''}.
\end{eqnarray}}
The last inequality holds by the facts $\log x>(\log\log x)^2$ for $x>e$,
$\log x>\log\log x$ for $x>1$ and by the inequalities \eqref{eq3.56}-\eqref{eq3.58}.
Noting that $\left |e^x-1 \right |<|x|e^{|x|}$ for $x\in\R$
and $e^x>1+x+\frac {1}{2}x^2$ for $x>0$, for $-d_n<x<0$, we have
\begin{eqnarray}
\label{eq3.72}
&&
\Lambda(x) \left | A_{n,\lambda(x)} - 1 \right |
\nonumber
\\
&<&
\Lambda(x) \left | h_{n,\lambda}(x) \right |
e^{\left | h_{n,\lambda}(x) \right |}
\nonumber
\\
&<&
\Lambda(x)\left[1+\left(1+\lambda^2\right)a_n^2x\right]^{-2}
e^{-x}\left(1+\lambda^2\right)a_n^2
\nonumber
\\
&&
\cdot
\left\{\frac {\left(1+\lambda^2\right)^3}{\lambda^2}
\left[1+\left(1+\lambda^2\right)a_n^2x\right]^{-2}+
\frac {1}{2}x^2+2 \left(1+\lambda^2\right)a_n^2x^2+2|x|\right\}
e^{\mathbb{C}_{5,\lambda}^{''}}
\nonumber
\\
&=&
a_n^2 \left[1+\left(1+\lambda^2\right)a_n^2x\right]^{-2}
\left(1+\lambda^2\right)
\exp\left\{-e^{-x}-x+\mathbb{C}_{5,\lambda}^{''}\right\}
\nonumber
\\
&&
\cdot
\left\{\frac {\left(1+\lambda^2\right)^3}{\lambda^2}
\left[1+\left(1+\lambda^2\right)a_n^2x\right]^{-2}+
\frac {1}{2}x^2+2\left(1+\lambda^2\right)a_n^2x^2+2|x|\right\}
\nonumber
\\
&<&
a_n^2 \left[1-\left(1+\lambda^2\right)a_n^2d_n\right]^{-2}
\left(1+\lambda^2\right)
\exp\left\{-\frac {1}{2}x^2+\mathbb{C}_{5,\lambda}^{''}-1\right\}
\nonumber
\\
&&
\cdot
\left\{\frac {\left(1+\lambda^2\right)^3}{\lambda^2}
\left[1 - \left(1+\lambda^2\right)a_n^2d_n\right]^{-2}+
\frac {1}{2}x^2+2\left(1+\lambda^2\right)a_n^2d_n^2+2|x|\right\}
\nonumber
\\
&=&
a_n^2 \left[1 - \left(1+\lambda^2\right)a_n^2d_n\right]^{-2}
\left(1+\lambda^2\right) e^{\mathbb{C}_{5,\lambda}^{''}-1}
\nonumber
\\
&&
\cdot
\left\{ \frac {\left(1+\lambda^2\right)^3}{\lambda^2}
\left[ 1 - \left(1+\lambda^2\right)a_n^2d_n\right]^{-2}e^{-\frac {1}{2}x^2}+
\frac {1}{2}x^2e^{-\frac {1}{2}x^2}+
2\left(1+\lambda^2\right)a_n^2d_n^2e^{-\frac {1}{2}x^2}+
2|x|e^{-\frac {1}{2}x^2}\right\}
\nonumber
\\
&<&
a_n^2 \left[1 - \left(1+\lambda^2\right)a_n^2d_n\right]^{-2}
\left(1+\lambda^2\right) e^{\mathbb{C}_{5,\lambda}^{''}-1}
\nonumber
\\
&&
\cdot
\left\{ \frac {\left(1+\lambda^2\right)^3}{\lambda^2}
\left[ 1 - \left(1+\lambda^2\right)a_n^2d_n\right]^{-2}+1+2\left(1+\lambda^2\right)
a_n^2\left[\log\left( \left(1+\lambda^2\right) b_n^2\right)\right]^2+1\right\}
\nonumber
\\
&<&
\mathbb{C}_{6,\lambda}^{''}a_n^2.
\end{eqnarray}
Thus putting \eqref{eq3.72} into \eqref{eq3.66} we can show that \eqref{eq3.53} holds.

Now the remainder is to show that \eqref{eq3.54} holds.
For $-\infty<x\leq-d_n$, noting that $e^x>1+x$ for $x \in \R$
and the values of $a_{n}$ and $d_{n}$, we have
\begin{eqnarray}
\label{eq3.73}
0\leq\Lambda(x)\leq\Lambda \left(-d_n\right) = \left(1+\lambda^2\right)a_n^2
\end{eqnarray}
and
\begin{eqnarray}
\label{eq3.74}
&&
e^{-\frac {1}{2}\left(1+\lambda^2\right)a_n^2d_n^2}
\left[1 - \left(1+\lambda^2\right)a_n^2d_n\right]^{-2}
\left\{1-\frac {\left(1+\lambda^2\right)^4}{\lambda^2}
a_n^2\left[ 1 - \left(1+\lambda^2\right)a_n^2d_n\right]^{-2}\right\}
\nonumber
\\
&>&
e^{-\frac {1}{2} \left(1+\lambda^2\right)a_n^2d_n^2}
\left\{1-\frac {\left(1+\lambda^2\right)^4}{\lambda^2}a_n^2
\left[ 1 - \left(1+\lambda^2\right)a_n^2d_n\right]^{-2}\right\}
\nonumber
\\
&>&
\left[1-\frac {1}{2} \left(1+\lambda^2\right)a_n^2d_n^2\right]
\left\{1-\frac {\left(1+\lambda^2\right)^4}{\lambda^2}
a_n^2\left[ 1 - \left(1+\lambda^2\right)a_n^2d_n\right]^{-2}\right\}
\nonumber
\\
&=&
1-\frac {\left(1+\lambda^2\right)^4}{\lambda^2}a_n^2
\left[ 1 - \left(1+\lambda^2\right)a_n^2d_n\right]^{-2}-
\frac {1}{2}\left(1+\lambda^2\right)a_n^2d_n^2+
\frac {\left(1+\lambda^2\right)^5}{2\lambda^2}a_n^4
d_n^2\left[ 1 - \left(1+\lambda^2\right)a_n^2d_n\right]^{-2}
\nonumber
\\
&>&
1-\frac {\left(1+\lambda^2\right)^4}{\lambda^2}a_n^2
\left[ 1 - \left(1+\lambda^2\right)a_n^2d_n\right]^{-2}-
\frac {1}{2} \left(1+\lambda^2 \right) a_n^2d_n^2.
\end{eqnarray}
Thus, by \eqref{eq1.5}, \eqref{eq1.6}, \eqref{eq3.56}, \eqref{eq3.57},
\eqref{eq3.58}, \eqref{eq3.74} and Lemma \ref{lem3.3},
for $-\infty<x\leq-d_n$, we have
\begin{eqnarray}
\label{eq3.75}
&&
F_\lambda^n \left(a_nx+b_n\right) \leq
F_\lambda^n \left(b_n-a_nd_n\right)
\nonumber
\\
&<&
\left\{ 1-\frac {2\phi \left(b_n-a_nd_n\right)
\phi\left(\lambda \left(b_n-a_nd_n\right)\right)}
{|\lambda| \left(1+\lambda^2\right)
\left(b_n-a_nd_n\right)^2}
\left[ 1-\frac {\left(1+\lambda^2\right)^2}{\lambda^2}
\left(b_n-a_nd_n\right)^{-2}\right] \right\}^n
\nonumber
\\
&=&
\left\{1-n^{-1}e^{d_n} e^{-\frac {1}{2}\left(1+\lambda^2\right)
a_n^2d_n^2}\left[1 - \left(1+\lambda^2\right)a_n^2d_n\right]^{-2}
\left[1-\frac {\left(1+\lambda^2\right)^4}{\lambda^2}
a_n^2\left[1 - \left(1+\lambda^2\right)a_n^2d_n\right]^{-2}\right] \right\}^n
\nonumber
\\
&<&
\exp\left\{-e^{d_n}e^{-\frac {1}{2}
\left(1+\lambda^2\right)a_n^2d_n^2}
\left[1 - \left(1+\lambda^2\right)a_n^2d_n\right]^{-2}
\left[1-\frac {\left(1+\lambda^2\right)^4}{\lambda^2}
a_n^2\left[ 1 - \left(1+\lambda^2\right)a_n^2d_n\right]^{-2}\right]\right\}
\nonumber
\\
&<&
\exp\left\{-e^{d_n}\left[1-\frac {\left(1+\lambda^2\right)^4}{\lambda^2}
a_n^2\left[ 1 - \left(1+\lambda^2\right)a_n^2d_n\right]^{-2} -
\frac {1}{2}\left(1+\lambda^2\right)a_n^2d_n^2\right]\right\}
\nonumber
\\
&=&
\exp\Bigg\{-\log\left[ \left(1+\lambda^2\right)b_n^2\right]+
\frac {\left(1+\lambda^2\right)^3}{\lambda^2}
\left[ 1 - \left(1+\lambda^2\right)a_n^2d_n\right]^{-2}
\left(1+\lambda^2\right)a_n^2
\log\left[ \left(1+\lambda^2\right) b_n^2\right]
\nonumber
\\
&&
+\frac {1}{2}\left(1+\lambda^2\right)
a_n^2\log\left[ \left(1+\lambda^2\right)b_n^2\right]
\left[\log\log\left( \left(1+\lambda^2\right) b_n^2\right)\right]^2\Bigg\}
\nonumber
\\
&<&
\exp\Bigg\{-\log\left[ \left(1+\lambda^2\right)b_n^2\right]+
\frac {\left(1+\lambda^2\right)^3}{\lambda^2}
\left[ 1 - \left(1+\lambda^2\right)a_n^2d_n\right]^{-2}
\left(1+\lambda^2\right)a_n^2
\log\left[ \left(1+\lambda^2\right)b_n^2\right]
\nonumber
\\
&&
+\frac {1}{2}\left(1+\lambda^2\right)
a_n^2\left[\log\left( \left(1+\lambda^2\right)b_n^2\right)\right]^2\Bigg\}
\nonumber
\\
&=&
a_n^2\left(1+\lambda^2\right)
\exp\Bigg\{\frac {\left(1+\lambda^2\right)^3}{\lambda^2}
\left[1 - \left(1+\lambda^2\right)a_n^2d_n\right]^{-2}
\left(1+\lambda^2\right)a_n^2\log\left[\left(1+\lambda^2\right)b_n^2\right]
\nonumber
\\
&&
+\frac {1}{2}\left(1+\lambda^2\right)
a_n^2\left[\log\left( \left(1+\lambda^2\right)b_n^2\right)\right]^2\Bigg\}
\nonumber
\\
&<&
\mathbb{C}_{7,\lambda}^{''}a_n^2.
\end{eqnarray}
Combining \eqref{eq3.73} with \eqref{eq3.75}, we  see that \eqref{eq3.54} holds.
The proof is complete.
\qed

\end{document}